\documentclass{article}
\usepackage{amssymb,latexsym}
\usepackage{graphicx}
\usepackage{amsmath}
\usepackage{amsthm}
\usepackage{amscd}
\usepackage{bm}
\begin{document}
\newcommand{\D}{\displaystyle}
\newcommand{\DF}[2]{\frac{\D#1}{\D#2}}
\def\Xint#1{\mathchoice
{\XXint\displaystyle\textstyle{#1}}%
{\XXint\textstyle\scriptstyle{#1}}%
{\XXint\scriptstyle\scriptscriptstyle{#1}}%
{\XXint\scriptscriptstyle\scriptscriptstyle{#1}}%
\!\int}
\def\XXint#1#2#3{{\setbox0=\hbox{$#1{#2#3}{\int}$}
\vcenter{\hbox{$#2#3$}}\kern-.5\wd0}}
\def\ddashint{\Xint=}
\def\dashint{\Xint-}
\pagestyle{plain}
\newpage
\normalsize
\title{A
regularity theory for multiple-valued Dirichlet minimizing maps
\thanks{Department of mathematics, Rice University, Houston, TX 77005; weizhu@math.rice.edu}}
\author{Wei Zhu}
\maketitle

\begin{abstract}
This paper discusses the regularity of multiple-valued Dirichlet
minimizing maps into the sphere. It shows that even at branched
point, as long as the normalized energy is small enough, we have
the energy decay estimate. Combined with the previous work by
Chun-Chi Lin, we get our first estimate that
$\mathcal{H}^{m-2}(\mbox{singular set})=0$. Furthermore, by
looking at the tangent map and using dimension reduction argument,
we show that the singular set is at least of codimension 3.
\end{abstract}
\newpage
\tableofcontents
\newpage
\newtheorem{theorem}{Theorem}[section]
\newtheorem{definition}{Definition}[section]
\newtheorem{remark}{Remark}[section]
\newtheorem{corollary}{Corollary}[section]
\newtheorem{proposition}{Proposition}[section]
\newtheorem{lemma}{Lemma}[section]
\section{Introduction}
The regularity of harmonic maps between Riemannian manifolds has
been a fascinating subject in recent years. The very first general
result on this is due to $\cite{su1}$, in which they proved that a
bounded, energy minimizing map $u:M^n\to N^k$ is regular (in the
interior) except for a closed set $S$ of Hausdorff dimension at
most $n-3$. One important technique they use in the paper is for
lowering the dimension of $S$ under the condition that certain
smooth harmonic maps of spheres into $N$ are trivial. This can be
checked in some interesting cases, for example if $N$ has
nonpositive curvature. They showed $S=\emptyset$, i.e, any energy
minimizing map into such a manifold is smooth. Use that method,
they $\cite{su2}$ are also able to reduce the dimension of $S$ if
$N$ is a sphere. The result is as follows:
\begin{theorem}[$\cite{su2}$, Theorem 2.7] For $k\ge 2$, define a
number $d(k)$ by setting
$$d(2)=2,d(3)=3$$
$$d(k)=[\min\{\DF{k}{2}+1,6\}]\;\mbox{for}\;k\ge 4$$
where $[\cdot]$ denotes the greatest integer in a number. If $n\le
d(k)$, then every energy minimizing map from a manifold $M$ of
dimension $n$ into $\mathbb{S}^k\subset\mathbb{R}^{k+1}$ is smooth
in the interior of $M$. If $n=d(k)+1$, such a map has at most
isolated singularities, and in general the singular set is a
closed set of Hausdorff dimension at most $n-d(k)-1$.
\end{theorem}
This same question in liquid crystal configurations
setting($n=3,k=2$) has been studied independently by
Hardt-Kinderlehrer-Lin using blowing-up argument in $\cite{hkl}$.
\\A few years later, Theorem 1.1 was extended to stable-stationary
harmonic maps $u\in H^1(\Omega,\mathbb{S}^k)$, $k\ge 3$ by
Hong-Wang $\cite{hw}$. Stable-stationary harmonic maps are
harmonic maps with zero domain first variation and nonnegative
range second variation. Examples of stable-stationary harmonic
maps include
energy minimizing maps.\\
In a recent work of Lin-Wang $\cite{lw}$, they improved the
theorems by $\cite{su2},\cite{hw}$ for $4\le k\le 7$ as follows:
\begin{theorem}[\cite{lw}, Theorem 1]
Define
\begin{equation*}
\tilde{d}(k)=
\begin{cases} 3 & k=3\\
4 & k=4\\
5 & 5\le k\le 9\\
6 & k\ge 10.
\end{cases}
\end{equation*}
For $k\ge 3$, let $u\in H^1(\Omega, \mathbb{S}^k)$ be a
stable-stationary harmonic map, then the singular set $S$ of $u$
has Hausdorff dimension at most $n-\tilde{d}(k)-1$.
\end{theorem}
We can also talk about the energy minimizing problems in the
setting of multiple-valued functions (maps) thanks to the
monumental work $\cite{af}$. After Almgren gave suitable
definitions of derivative and Sobolev space for multiple-valued
functions, the question of minimizing energy among functions with
the prescribed boundary data becomes legitimate. Furthermore, he
was able to show that any Dirichlet energy minimizing
multiple-valued function is regular in the interior and has branch
point of codimension at least 2. Although the primary purpose in
$\cite{af}$ of introducing multiple-valued functions is to
approximate almost flat mass-minimizing integral currents by
graphs of Dirichlet minimizing multiple-valued functions, the
subject of multiple-valued functions in the sense of Almgren turns
out to be also interesting in its own. See some recent works
$\cite{cs},\cite{gj},\cite{lc1},\cite{lc2},\cite{mp},\cite{zw1},\cite{zw2},\cite{zw3}$.\\
In the same spirit of $\cite{af}$, Chun-Chi Lin (in $\cite{lc1}$)
considered the energy minimizing multiple-valued map into spheres.
Specifically, he showed that for points not in the branch set
$B_0$, as long as the normalized energy is small, the map is
regular there (see more of this discussion in section three). We
will continue his work by examining the local behavior of points
in $B_0$. The main idea is to use the blowing-up analysis at this
point. The blowing-up sequence converges strongly to a Dirichlet
minimizing function which is regular due to $\cite{af}$. Hence it
guarantees the energy of the original map near this point
satisfies some growth condition. Combining our result with the
result in $\cite{lc1}$, we conclude that the minimizing map is
regular at any point as long as the normalized energy there is
small enough thanks to Morrey's growth lemma. This gives us the
first $m-2$ estimate. Then, using dimension reduction argument, we
get our main result:
\begin{theorem}
Let
$u\in\mathcal{Y}_2(\mathbb{B}_1^m(0),\underline{Q}(\mathbb{S}^{n-1}))$
be a strictly defined, Dirichlet minimizing map. Then it is
H$\ddot{o}$lder continuous away from the boundary except for a
closed subset $S\subset \mathbb{B}_1^m(0)$ such that $\dim(S)\le
m-3$.
\end{theorem}
The assumption that we are looking at points in $B_0$ is important
in the blowing-up process because we need to get suitable constant
of the form $Q[[b]]$ for some $b\in \mathbb{S}^{n-1}$ in order for
the
subtraction between two Q-tuples to make sense.\\
There are some other interesting questions which are not addressed
in this paper, and still open to the author's knowledge. A first
one is whether our result is an optimal one. We are hoping to have
some similar results as in $\cite{su2},\cite{lw}$. Some new
techniques are expected because $\cite{su2} \cite{lw}$ both use
Bochner formula, which is no longer available in the
multiple-valued functions setting. \\
A second one is the regularity for stationary harmonic
multiple-valued functions. There are already some positive results
for
this in the two dimensional case, see $\cite{lc2}$.\\
Another one is the branching behavior. Chun-Chi Lin (in
$\cite{lc1}$) has done some work on this. But there was some
problem with that. Basically speaking, the monotonicity formula
for frequency function he used in his proof actually does not
necessarily hold for multiple-valued maps. Some new idea is
probably needed to get around this
obstacle.\\
It is my great pleasure to thank my thesis advisor Professor
Robert Hardt for his support, encouragement and kindness during
the years at Rice. A lot of this work was stimulated by
$\cite{hkl}$.
\section{Preliminaries}
Most of the notations, definitions and known results about
multiple-valued functions that we need can be found in
$\cite{zw1}$. The reader is also referred to $\cite{af}$ for more
details. We also use standard terminology in geometric measure
theory,
all of which can be found on page 669-671 of the treatise {\it Geometric Measure Theory} by H. Federer $\cite{fh}$.\\
For reader's convenience, here we state some useful results not
included in $\cite{zw1}$. The proofs of them can be found in
$\cite{af}$.
\begin{theorem}[$\cite{af}, \S 2.6$]
(a) $0<r_0<\infty$.\\
(b) $A\subset \mathbb{R}^m$ is connected, open, and bounded with $\mathbb{B}_{r_0}^m(0)\subset A$. $\partial A$
is an $m-1$ dimensional submanifold of $\mathbb{R}^m$ of class 1.\\
(c) $f:A\rightarrow \mathbb{Q}$ is strictly defined and is Dir minimizing.\\
(d) $D,H,N:(0,r_0)\rightarrow\mathbb{R}$ are defined for $0<r<r_0$ by setting
\begin{equation*}
\begin{split}
&D(r)=\mbox{Dir}(f;\mathbb{B}_r^m(0))\\
&H(r)=\int_{\partial \mathbb{B}_r^m(0)}\mathcal{G}(f(x),Q[[0]])^2
d\mathcal{H}^{m-1} x\\
&N(r)=rD(r)/H(r)\;\mbox{provided}\;H(r)>0.
\end{split}
\end{equation*}
(e) $\mathcal{N}: A\to \mathbb{R}$ is defined for $x\in A$ by
setting
$$\mathcal{N}(x)=\lim_{r\downarrow 0}
r\mbox{Dir}(f;\mathbb{B}_r^m(x))/\int_{\partial
\mathbb{B}_r^m(x)}\mathcal{G}(f(z),Q[[0]])^2d\mathcal{H}^{m-1}z$$
provided this limit exists.\\ (f) $H(r)>0$ for some $0<r<r_0$.\\
{\bfseries Conclusions.}\\
For $\mathcal{L}^1$ almost all $0<r<r_0$,\\
Squash Deformation:\\
$$D(r)=\int_{x\in\partial \mathbb{B}_r^m(0)}<\xi_0\circ f(x),D(\xi_0\circ f)(x,x/|x|)>d\mathcal{H}^{m-1}x$$
Squeeze Deformation:\\
$$r\cdot D'(r)=(m-2)D(r)+2r\int_{x\in\partial\mathbb{B}_r^m(0)}|D(\xi_0\circ f)(x,x/|x|)|^2 d\mathcal{H}^{m-1} x$$
\end{theorem}
\begin{remark}(1) For convenience, we will use $\partial f/\partial r$ to represent $Df(x,x/|x|)$ for any multiple-valued function $f$ whenever
the derivative exits.\\
(2) Noticing that the squeeze deformation comes from a domain deformation, the squeeze deformation formula still holds for multiple-valued maps.\\
(3) We can replace those $\xi_0$ by $\xi$ because that in the
proof of those formulas, the only things we need for $\xi_0$ are
that $\xi_0(\lambda x)=\lambda \xi_0(x)$, for $\lambda>0$,
$x\in\mathbb{Q}$, and $D(r)=\int_{B_r^m(0)}|D(\xi_0\circ f)|^2 d
\mathcal{L}^m$, both of which still hold for $\xi$.
\end{remark}
\begin{theorem}[$\cite{af}, \S 2.9$]
Corresponding to numbers $0<s_0<\infty$, $1<K<\infty$, and(not
necessarily distinct points)
$q_1,\cdot\cdot\cdot,q_Q\in\mathbb{R}^n$ we can find
$J\in\{1,\cdot\cdot\cdot,Q\},k_1,\cdot\cdot\cdot,k_J\in\{1,\cdot\cdot\cdot,Q\}$,
distinct
points $p_1,\cdot\cdot\cdot,p_J\in\{q_1,\cdot\cdot\cdot,q_Q\}$, and $s_0\le r\le Cs_0$ such that\\
(1) $|p_i-p_j|>2Kr$ for each $1\le i<j\le J$,\\
(2) $\mathcal{G}(\sum_{i=1}^Q [[q_i]],\sum_{i=1}^J k_i[[p_i]])\le Cs_0/(Q-1)^{1/2},$\\
(3) $z\in\mathbb{Q}$ with $\mathcal{G}(z,\sum_{i=1}^Q[[q_i]])\le
s_0$ implies $\mathcal{G}(z,\sum_{i=1}^J k_i[[p_i]])
\le r$,\\
(4) in case $J=1$, diam(spt($\sum_{i=1}^Q [[q_i]]))\le
Cs_0/(Q-1)$; here
$$C=1+[(2K)(Q-1)^2]^1+[(2K)(Q-1)^2]^2+\cdot\cdot\cdot+[(2K)(Q-1)^2]^{Q-1}.$$
\end{theorem}
\begin{theorem}[$\cite{af}, \S 2.10$]
Corresponding to\\
(a) $J\in\{1,2,\cdot\cdot\cdot,Q\}$,\\
(b) $k_1,k_2,\cdot\cdot\cdot,k_J\in\{1,2,\cdot\cdot\cdot,Q\}$ with
$k_1+k_2+\cdot\cdot\cdot+k_J=Q$,\\
(c) distinct points
$p_1,p_2,\cdot\cdot\cdot,p_J\in\mathbb{R}^n$,\\
(d) $0<s_1<s_2=2^{-1}\inf\{|p_i-p_j|:1\le i<j\le J\}$,\\ we set
$$\mathbb{P}=\mathbb{Q}\cap\{\sum_{i=1}^Q[[q_i]]:q_1,\cdot\cdot\cdot,q_Q\in\mathbb{R}^n\;\mbox
{with card}\{i:q_i\in\mathbb{B}_{s_1}^n(p_j)\}=k_j\;\mbox{for
each}\;j=1,\cdot\cdot\cdot,J\}$$ {\bfseries Conclusions}:\\ There
is a map $\Phi:\mathbb{Q}\to\mathbb{P}$ such that\\
(1) $\Phi(q)=q$ whenever $q\in\mathbb{Q}$ with
$\mathcal{G}(q,\sum_{i=1}^J k_i[[p_i]])\le s_1$,\\
(2) $\Phi(q)=\sum_{j=1}^J k_j[[p_j]]$ whenever $q\in\mathbb{Q}$
with $\mathcal{G}(q,\sum_{i=1}^J k_i[[p_i]])\ge s_2$,\\
(3) $\mathcal{G}(q,\Phi(q))\le \mathcal{G}(q,\sum_{i=1}^J
k_i[[p_i]])$ for each $q\in\mathbb{Q}$,\\
(4) Lip $\Phi\le 1+Q^{1/2}s_1/(s_2-s_1).$
\end{theorem}
\begin{theorem}[$\cite{af}, \S 2.13$]
(a) In case $m=2$, $\omega_{2.13}=1/Q$.\\
(b) In case $m\ge 3, 0<\epsilon_Q<1$ is as defined as in $\cite{af},\S 2.11$ and $0<\omega_{2.13}<1$ is defined by the requirement
$$m-2+2\omega_{2.13}=(m-2)(1+\epsilon_Q)/(1-\epsilon_Q).$$
(c) $$\Gamma_{2.13}=4^{1-\omega_{2.13}}[2^{m/2}/(1-2^{-\omega_{2.13}})+3\cdot 2^{m-1+\omega_{2.13}}](m\alpha(m))^{-1/2}
\mbox{Lip}(\xi)^2\mbox{Lip}(\xi^{-1}).$$
(d) $f\in\mathcal{Y}_2(\mathbb{R}^m,\mathbb{Q})$ is strictly defined and $f|\mathbb{B}_1^m(0)$ is Dir minimizing with
Dir$(f;\mathbb{B}_1^m(0))>0$.\\
{\bfseries Conclusions.}\\
(1) For each $z\in \mathbb{B}_1^m(0),0<r<1-|z|$, and $0<s\le 1$,
$$\mbox{Dir}(f;\mathbb{B}_{sr}^m(z))\le s^{m-2+2\omega_{2.13}}\mbox{Dir}(f;\mathbb{B}_r^m(z)).$$
(2) Whenever $0<\delta<1$ and $p,q\in\mathbb{B}_{1-\delta}^m(0)$,
$$\mathcal{G}(f(p),f(q))\le \Gamma_{2.13}\delta^{-m/2}\mbox{Dir}(f;\mathbb{B}_1^m(0))^{1/2}|p-q|^{\omega_{2.13}},$$
in particular, $f|\mathbb{B}_{1-\delta}^m(0)$ is H$\ddot{o}$lder continuous with exponent $\omega_{2.13}$. \\
(3) Corresponding to each bounded open set $A$ such that $\partial A$ is a compact $m-1$ dimensional submanifold of $\mathbb{R}^m$ of
class 1, there is a constant $0<\Gamma_A<\infty$ with the following property. Whenever $g\in\mathcal{Y}_2(A,\mathbb{Q})$ is Dir
minimizing and $p,q\in A$,
$$\mathcal{G}(g(p),g(q))\le \Gamma_A\mbox{Dir}(g;A)^{1/2}\sup\{\mbox{dist}(p,\partial A)^{-m/2},\mbox{dist}(q,\partial A)^{-m/2}\}
|p-q|^{\omega_{2.13}}.$$
\end{theorem}
\section{Some Remarks on $\cite{lc1}$}
In $\cite{lc1}$, Chun-Chi Lin introduced the set\\
$B_0=\{x\in \mathbb{B}_1^m(0):\mbox{a Lebesgue point of}\;f,
\xi^{-1}\circ \rho\circ
AV_{r,x}(\xi\circ f)=Q[[b_r]], $\\
$\mbox{for any small enough radius}\; r>0, b_r\in \mathbb{R}^n\;
\mbox{and}\;
AV_{r,x}(\xi\circ f)=\dashint_{\partial \mathbb{B}_r^m(x)} \xi\circ f \}.$\\
He proved that for a point not in $B_0$, if the normalized energy
of $f$ is small enough there, then the energy of $f$ near this
point satisfies some growth condition. The key ingredients are the
induction on $Q$ and finding a comparison map. In order to use the
induction, we need $J\ge 2$. This is
guaranteed by our assumption that the point we are looking at is not in $B_0$. He did not explain that in his paper. Here is the detail:\\
If $a\not\in B_0$, i.e. there is $r>0$, such that $\xi^{-1}\circ \rho\circ AV_{r,a}(\xi\circ f)\not=Q[[b]]$ for any $b\in \mathbb{R}^n$.
We may as well just assume that
$$\xi^{-1}\circ\rho\circ AV_{1,a}(\xi\circ f)\not= Q[[b]]\;\mbox{for any}\;b\in \mathbb{R}^n.$$
Now instead of letting $q^*\in \mathbb{Q}^*$ be the point in
$\mathbb{Q}^*$ such that
$$|q^*-AV_{1,a}(\xi\circ f)|=\mbox{dist}(AV_{1,a}(\xi\circ f),\mathbb{Q}^*)$$
we let $q^*=\rho\circ AV_{1,a}(\xi\circ f)$, and $q=\sum_{i=1}^Q
[[q_i]]=\xi^{-1}(q^*)$.  With these points
$q_1,q_2,\cdot\cdot\cdot,q_Q$, $1<K<\infty$ and constant $s_0$ to
be chosen later, we find $J\in\{1,2,\cdot\cdot\cdot,Q\}$,
$k_1,k_2,\cdot\cdot\cdot,k_J\in\{1,2,\cdot\cdot\cdot,Q\}$,
distinct points $p_1,\cdot\cdot\cdot,p_J\in \mathbb{R}^n$ as in Theorem 2.2. Let $q_0=\sum_{i=1}^J k_i[[p_i]]$.\\
If $J=1$, from Theorem 2.2 (4)
\begin{equation*}
\begin{split}
\mbox{diam(spt}(\sum_{i=1}^Q[[q_i]]))&=\mbox{diam(spt}(\xi^{-1}\circ\rho\circ
AV_{1,a}(\xi\circ f))) \\
&\le Cs_0/(Q-1); \end{split}
\end{equation*} but we already know that $\xi^{-1}\circ\rho\circ
AV_{1,a}(\xi\circ f)\not=Q[[b]]$ for any $b\in \mathbb{R}^n$,
hence diam(spt($\sum_{i=1}^Q[[q_i]]$)) is a fixed positive number.
So we can choose
$s_0$ small enough to guarantee that $J\ge 2$.\\
We also have to show that the rest of the proof in $\cite{lc1}$ is
still valid after we choose the different $q^*$. This is because
the only place where $q^*$ is used in $\cite{lc1}$ is to show
$$\int_{\partial \mathbb{B}_1^m(a)}|\xi\circ f(x)-q^*|^2\le C\int_{\partial
\mathbb{B}_1^m(a)} |\xi\circ f(x)-AV_{1,a}(\xi\circ
f)|^2\;\mbox{for some constant}\;C.$$
We still have this because\\
$\int_{\partial \mathbb{B}_1^m(a)} |\xi\circ f-q^*|^2=\int_{\partial \mathbb{B}_1^m(a)}|\xi\circ f(x)-\rho\circ AV_{1,a}(\xi\circ f)|^2$\\
$=\int_{\partial \mathbb{B}_1^m(a)}|\rho\circ\xi\circ
f(x)-\rho\circ AV_{1,a}(\xi\circ
f)|^2\le\;(\mbox{Lip}\;\rho)^2\int_{\partial \mathbb{B}_1^m(a)}
|\xi\circ f(x)-AV_{1,a}(\xi\circ f)|^2.$\\
Another thing that worth mentioning is in the proof of Lemma 4 in
$\cite{lc1}$, more precisely (2.12). He was claiming that $g_j$ is
H$\ddot{o}$lder continuous hence having growth condition on the
energy. But in fact since his work is only on points outside
$B_0$, and we do not know whether the origin is inside or outside
of the set $B_0$ for each $g_j$, the induction seems to be a
problem. However, using our result on branched points, we can
overcome this. Let's look at our result Theorem 7.1 in advance
(notice that our proof does not depend on induction or the result
in $\cite{lc1}$) , which says that at branched point,
$$D(r)\le Cr^{m-2+\omega_{2.13}},$$
for some constant $C$ depending on the dimensions and the total
energy $D(1)$. Now we claim that for each $g_j$ in $\cite{lc1}$,
there exists a positive constant $\alpha$ such that
$$D_g(r(1-t_Q))=\sum_{j=1}^J D_{g_j}(r(1-t_Q))\le
C(\alpha,m,n,Q,\mbox{total energy of f})r^{m-2+\alpha}.$$ This is
because if the origin is not in the corresponding set $B_0$ of
$g_j$, then the induction argument gives us the above estimate.
Otherwise, our result applies.\\
Finally, we modify the end of proof of Lemma 4 in $\cite{lc1}$ as
following: (reason that the original proof did not work is that by
considering two cases, the integration did not necessarily work)\\
Now we have
$$D_f(r)\le \DF{8}{7}Cr^{m-2+\alpha}+\DF{1}{7(m-1)}rD'_f(r).$$
Let's denote $D_f(r)$ by $\phi(r)$. The original inequality
becomes
$$\phi(r)'-\DF{\phi(r)}{r}7(m-1)+8Cr^{m-3+\alpha}(m-1)\ge 0.$$
Multiply both sides by $r^{-7(m-1)}$, we get
$$\DF{d}{dr}[\phi(r)r^{-7(m-1)}+\DF{8C(m-1)}{5-6m+\alpha}r^{5-6m+\alpha}]\ge 0.$$
Hence
$$\phi(r)r^{-7(m-1)}+\DF{8C(m-1)}{5-6m+\alpha}r^{5-6m+\alpha}\le
\phi(1)+\DF{8C(m-1)}{5-6m+\alpha}:=M,$$
\begin{equation*}
\begin{split}
D_f(r)=\phi(r)&\le
(M+\DF{8C(m-1)}{6m-5-\alpha}r^{5-6m+\alpha})r^{7(m-1)}\\
&=Mr^{7(m-1)}+\DF{8C(m-1)}{6m-5-\alpha}r^{m-2+\alpha}\\
&\le \max(M,\DF{8C(m-1)}{6m-5-\alpha})r^{m-2+\alpha}
\end{split}
\end{equation*}
while the last inequality follows because $7(m-1)>m-2+\alpha$.
\section{Maximum Principle for Multiple-Valued Dirichlet Minimizing Functions}
\begin{lemma} Given a positive number $M$, and $\epsilon>0$, define the
retraction function $\Pi_M$ as follows
\begin{equation*}
\Pi_M(x)=
\begin{cases} x & |x|\le M\\
\DF{x}{|x|}M & \mbox{otherwise}
\end{cases}
\end{equation*}
Let $T=\{x:|x|\ge M+\epsilon\}$. Then Lip($\Pi_M|T)\le \DF{M}{M+\epsilon}$.
\end{lemma}
\begin{proof} For $x,y\in T$,
\begin{equation*}
\begin{split}
|\Pi_M(x)-\Pi_M(y)|^2 & =|\Pi_M(x)|^2+|\Pi_M(y)|^2-2<\Pi_M(x),\Pi_M(y)>\\
&=2M^2-2M^2\DF{<x,y>}{|x||y|}\\
&=\DF{M^2}{|x||y|}(2|x||y|-2<x,y>)\\
&\le\DF{M^2}{|x||y|}(|x|^2+|y|^2-2<x,y>)\\
&=\DF{M^2}{|x||y|}|x-y|^2\le\DF{M^2}{(M+\epsilon)^2}|x-y|^2.
\end{split}
\end{equation*}
\end{proof}
\begin{definition} For a Q-valued function $f$, define
$$|f(x)|:=\mbox{Max}\{|f_1(x)|,|f_2(x)|,\cdot\cdot\cdot,|f_Q(x)|\},$$
where $f=\sum_{i=1}^Q [[f_i]]$.
\end{definition}
\begin{theorem}
If $f:\mathbb{B}_1^m(0)\rightarrow \mathbb{Q}$ is strictly defined and Dir minimizing with boundary data
$g:\partial\mathbb{B}_1^m(0)\rightarrow \mathbb{Q}$, where $f\in \mathcal{Y}_2(\mathbb{B}_1^m(0),\mathbb{Q}),g\in
\partial\mathcal{Y}_2(\partial \mathbb{B}_1^m(0),\mathbb{Q})$, then
$$\sup_{x\in \mathbb{B}_1^m(0)} |f(x)|\le  \sup_{x\in \partial\mathbb{B}_1^m(0)}|g(x)|.$$
\end{theorem}
\begin{proof}
We may assume that $M:=\sup_{x\in\partial\mathbb{B}_1^m(0)}|g(x)|<
\infty$. \\
If the statement is not true, i.e. there is a point $x_0\in \mathbb{B}_1^m(0)$, such that
$|f(x_0)|>M$.
\newline
Claim: $f(x)=f(x_0)$ for all $x\in \mathbb{B}_1^m(0)$.\\
Define the set $S=\{x\in\mathbb{B}_1^m(0):f(x)=f(x_0)\}$, which is
not empty by the assumption. Since $f$ is continuous from Theorem
2.4, $S$ is closed in $\mathbb{B}_1^m(0)$.\\
Let $\Pi_M$ be the retraction function from $\mathbb{R}^m$ to
$\mathbb{B}_M^m(0)$. Consider the new comparison Q-valued function
$h=(\Pi_M)_{\sharp}\circ f$, which has boundary
data $g$ and whose energy is no more than that of $f$ because Lip$(\Pi_M)\le 1$.\\
Take any point $y\in S$, because of the continuity of $f$, there
is a neighborhood $U$ of $y$ in $\mathbb{B}_1^m(0)$ such that
$|f(x)|\ge M+\epsilon, x\in U$, for some
$\epsilon$ small enough. \\
From Lemma 4.1, we know that Lip$(\Pi_M|U)\le \DF{M}{M+\epsilon}$,
hence
$$\mbox{Dir}(h;U)\le \DF{M}{M+\epsilon}\;\mbox{Dir}(f;U).$$
Therefore $f$ must be constant in $U$(otherwise, its energy is
nonzero. But the energy of $h$ in $U$ is strictly smaller than
that of $f$, contradicting to the fact that $f$ is Dir
minimizing). So $S$ is also open in $\mathbb{B}_1^m(0)$.
Therefore, $S=\mathbb{B}_1^m(0)$, which is a contradiction to the
assumption that $f|\partial \mathbb{B}_1^m(0)=g$.
\end{proof}
\section{Hybrid Inequality}
From now on, $m\ge 2$ and $n\ge 2$.
\begin{lemma} If $u:\mathbb{B}^m_1(0)\rightarrow
\underline{Q}(\mathbb{S}^{n-1})$ is strictly defined and Dir
minimizing, then for $a.e. \;0<r<1$,
$$\int_{\partial\mathbb{B}_r^m(0)} |\DF{\partial u}{\partial r}|^2 dx\le \int_{\partial\mathbb{B}_r^m(0)} |\nabla_{\mbox{tan}} u|^2 dx$$
\end{lemma}
\begin{proof}
For minimizing maps, we still have the squeeze formula:
$$r\cdot D'(r)=(m-2)\cdot D(r)+2r\cdot \int_{\partial \mathbb{B}_r^m(0)} |\DF{\partial u}{\partial r}|^2 d\mathcal{H}^{m-1}.$$
Noticing that $D'(r)=\int_{\partial \mathbb{B}_r^m(0)} |\DF{\partial u}{\partial r}|^2 d\mathcal{H}^{m-1}+\int_{\partial
\mathbb{B}_r^m(0)}|\nabla_{
\mbox{tan}}u|^2
d\mathcal{H}^{m-1}$, we have
$$2r\cdot \mbox{dir}(f,\partial \mathbb{B}_r^m(0))=(m-2)D(r)+r\cdot D'(r)$$
Therefore, $r\cdot D'(r)\le 2r\cdot \mbox{dir}(f,\partial \mathbb{B}_r^m(0))$, i.e. $D'(r)\le 2\mbox{dir}(f,\partial
\mathbb{B}_r^m(0))$.\\
Writing that in integration form,
$$\int_{\partial \mathbb{B}_r^m(0)}|\nabla_{\mbox{tan}} u|^2+\int_{\partial \mathbb{B}_r^m(0)}|\DF{\partial u}{\partial r}|^2\le
2\int_{\partial \mathbb{B}_r^m(0)}|\nabla_{\mbox{tan}} u|^2,$$
gives us the desired inequality.
\end{proof}
\begin{theorem}[Hybrid Inequality]
There is a positive constant $C$, depending only on $m,n,Q$ such that if $0<\lambda<1$, and if $u$ is a Dir-minimizer
in $\mathcal{Y}_2(\mathbb{B}_1^m(0),\underline{Q}(S^{n-1}))$, then
$$E_{1/2}(u)\le \lambda E_1(u)+C\lambda^{-1}\int_{\mathbb{B}_1^m(0)}|\xi\circ u-\overline{\xi\circ u}|^2 dx,$$
where $E_r(u)=r^{2-m}\int_{\mathbb{B}_r^m(0)} |Du|^2 dx$, $\overline{\xi\circ u}=\dashint_{\mathbb{B}_1^m(0)} (\xi\circ u )dx$.
\end{theorem}
\begin{proof}
For an increasing function $\eta$ on $[0,1]$,
$$\{s:\eta'(s)\ge 8[\eta(1)-\eta(0)]\}$$
has Lebesgue measure $\le 1/8$. In particular, there is an $r\in [1/2,1]$ such that $u|\partial \mathbb{B}_r^m(0)\in \partial
\mathcal{Y}_2(\partial \mathbb{B}_r^m(0),\underline{Q}
(\mathbb{S}^{n-1}))$,
$$\int_{\partial \mathbb{B}_r^m(0)}|\nabla_{\mbox{tan}} u|^2 d\mathcal{H}^{m-1}\le 8\int_{\mathbb{B}_1^m(0)}|Du|^2 d\mathcal{H}^m,$$
$$\int_{\partial \mathbb{B}_r^m(0)}|\xi\circ u-\overline{\xi\circ u}|^2 d\mathcal{H}^{m-1}\le 8 \int_{
\mathbb{B}_1^m(0)}|\xi\circ u-\overline{\xi\circ u}|^2 d
\mathcal{H}^m.$$
We claim that there exists a map $w\in \mathcal{Y}_2(\mathbb{B}_r^m(0),\underline{Q}(
\mathbb{S}^{n-1}))$ such that $w|\partial \mathbb{B}_r^m(0)=u|\partial \mathbb{B}_r^m(0)$ and
$$\int_{\mathbb{B}_r^m(0)}|Dw|^2 dx\le K(\int_{\partial \mathbb{B}_r^m(0)}|\nabla_{\mbox{tan}} u|^2 d\mathcal{H}^{m-1})^{1/2}(\int_{\partial
\mathbb{B}_r^m(0)}|
\xi\circ u-\overline{\xi\circ u}|^2 d\mathcal{H}^{m-1})^{1/2}$$
for some universal constant $K$.
This is because, we choose $h: \mathbb{B}_r^m(0)\rightarrow \mathbb{Q}$ such that $h|\partial
\mathbb{B}_r^m(0)=u|\partial \mathbb{B}_r^m(0)$ and $h$ is Dir-minimizing.
\begin{equation*}
\begin{split}
\int_{\mathbb{B}_r^m(0)}|Dh|^2 &=\int_{\partial \mathbb{B}_r^m(0)} <\xi\circ h(x),D(\xi\circ h)(x,\DF{x}{|x|})>d\mathcal{H}^{m-1}\\
&=\int_{\partial \mathbb{B}_r^m(0)}<\xi\circ h(x)-\overline{\xi\circ u},D(\xi\circ h)(x,\DF{x}{|x|})>d\mathcal{H}^{m-1}\\
&\le [\int_{\partial \mathbb{B}_r^m(0)} |\xi\circ h(x)-\overline{\xi\circ u}|^2 d\mathcal{H}^{m-1}]^{1/2}[\int_{\partial
\mathbb{B}_r^m(0)} |\DF{\partial h}{\partial r}|^2 d
\mathcal{H}^{
m-1}]^{1/2}\\
&=[\int_{\partial \mathbb{B}_r^m(0)} |\xi\circ u(x)-\overline{\xi\circ u}|^2 d\mathcal{H}^{m-1}]^{1/2}[\int_{\partial
\mathbb{B}_r^m(0)} |\DF{\partial h}{\partial r}|^2 d
\mathcal{H}^{
m-1}]^{1/2}
\end{split}
\end{equation*}
By Lemma 5.1, we have
$$\int_{\partial \mathbb{B}_r^m(0)}|\DF{\partial h}{\partial r}|^2\le \int_{\partial
\mathbb{B}_r^m(0)}|\nabla_{\mbox{tan}} h|^2=\int_{\partial \mathbb{B}_r^m(0)} |\nabla_{
\mbox{tan}} u|^2$$
Therefore
$$\int_{\mathbb{B}_r^m(0)}|Dh|^2\le [\int_{\partial \mathbb{B}_r^m(0)}|\xi\circ u-\overline{\xi\circ u}|^2d
\mathcal{H}^{m-1}]^{1/2}[\int_{\partial \mathbb{B}_r^m(0)}|
\nabla_{\mbox{tan}} u|^2d\mathcal{H}^{m-1}]^{1/2}$$
Unfortunately, the image of $h$ may not lie in $\underline{Q}(\mathbb{S}^{n-1})$. To correct this, we consider, for $a\in
\mathbb{B}_{1/2}^m(0)$, the projection
$$\Pi_a(x)=(x-a)/|x-a|,$$
and note that by Sard's Theorem, the composition
$(\Pi_a)_\sharp\circ h\in \mathcal{Y}_2(
\mathbb{B}_1^m(0),\underline{Q}(\mathbb{S}^{n-1}))$ for almost all a.\\
Using Fubini's Theorem, we estimate
$$\int_{
\mathbb{B}_{1/2}^m(0)}\int_{\mathbb{B}_r^m(0)}
|D((\Pi_a)_\sharp\circ h)(x)|^2 dxda$$ $$\le 4\int_{
\mathbb{B}_r^m(0)}|Dh(x)|^2\int_{\mathbb{B}_{1/2}^m(0)}(\mathcal{G}(h(x),Q[[a]]))^{-2}dadx$$
Now we claim that $\int_{\mathbb{B}_{1/2}^m(0)}(\mathcal{G}(h(x),Q[[a]]))^{-2} da<K$ for some universal constant $K$ independent of $x$.\\
Let $h(x)=\sum_{i=1}^Q [[h_i(x)]]$, then $(\mathcal{G}(h(x),Q[[a]]))^2=\sum_{i=1}^Q |h_i(x)-a|^2$. Hence
$$(\mathcal{G}(h(x),Q[[a]]))^{-2}\le |h_i(x)-a|^{-2},\;\mbox{for any}\;i$$
Applying Theorem 4.1 to the function $h$, and noticing that
$$h|\partial \mathbb{B}_r^m(0)\in
\partial\mathcal{Y}_2(\partial \mathbb{B}_r^m(0),\underline{Q}(\mathbb{S}^{n-1})),$$
we get
$$|h_i(x)|\le 1,\;\mbox{for any}\;i,\;\mbox{any}\;x\in \mathbb{B}_r^m(0)$$
Therefore
$$\int_{\mathbb{B}_{1/2}^m(0)}(\mathcal{G}(h(x),Q[[a]]))^{-2} da\le \int_{
\mathbb{B}_{1/2}^m(0)} |h_i(x)-a|^{-2} da$$
\begin{equation*}
\begin{split}
&=\int_{\mathbb{B}_{1/2}^m(-h_i(x))}|y|^{-2} dy,\;\mbox{by changing of variable}\;a=h_i(x)+y\\
&\le \int_{\mathbb{B}_2^m(0)}|y|^{-2} dy<\infty.
\end{split}
\end{equation*}
Hence $\int_{\mathbb{B}_{1/2}^m(0)}\int_{
\mathbb{B}_r^m(0)}|D((\Pi_a)_\sharp\circ h)(x)|^2 dxda\le K\int_{\mathbb{B}_r^m(0)}|Dh(x)|^2 dx$ for some constant $K$.
We may choose $a\in \mathbb{B}_{1/2}^m(0)$ such that $\int_{
\mathbb{B}_r^m(0)}|D((\Pi_a)_\sharp\circ h)(x)|^2 dx\le K\int_{\mathbb{B}_r^m(0)} |Dh(x)|^2 dx$. Letting
$w=[(\Pi_a|\mathbb{S}^{n-1})^{-1}]_\sharp\circ (\Pi_a)_\sharp\circ h$, we conclude that $w|\partial
\mathbb{B}_r^m(0)=u|\partial \mathbb{B}_r^m(0)$, and that
\begin{equation*}
\begin{split}
\int_{
\mathbb{B}_r^m(0)}|Dw(x)|^2 dx &\le [\mbox{Lip}(\Pi_a|\mathbb{S}^{n-1})^{-1}]^2\int_{
\mathbb{B}_r^m(0)}|D((\Pi_a)_\sharp\circ h)|^2 dx\\
&\le K\int_{\mathbb{B}_r^m(0)} |Dh|^2 dx
\end{split}
\end{equation*}
Now back to our desired result,
\begin{equation*}
\begin{split}
E_{1/2}(u) &=(1/2)^{2-m}\int_{
\mathbb{B}_{1/2}^m(0)}|Du|^2 dx\\
&\le 2^{m-2}\int_{\mathbb{B}_r^m(0)} |Du|^2 dx\\
&\le 2^{m-2}\int_{
\mathbb{B}_r^m(0)} |Dw|^2 dx\\
&\le 2^{m-2}K(\int_{\partial \mathbb{B}_r^m(0)}|\nabla_{\mbox{tan}} u|^2 d\mathcal{H}^{m-1})^{1/2}(\int_{\partial
\mathbb{B}_r^m(0)}|
\xi\circ u-\overline{\xi\circ u}|^2 d\mathcal{H}^{m-1})^{1/2}
\end{split}
\end{equation*}
Applying the inequality $ab\le \DF{1}{2}\delta a^2+\DF{1}{2}\delta^{-1}b^2$, with $\delta=\DF{\lambda}{2^mK}$, we have
\begin{equation*}
\begin{split}
E_{1/2}(u)&\le 2^{m-2}K(\DF{1}{2}\delta \int_{\partial
\mathbb{B}_r^m(0)}|\nabla_{\mbox{tan}} u|^2+\DF{1}{2}\delta^{-1}
\int_{\partial \mathbb{B}_r^m(0)} |\xi\circ u-\overline{\xi\circ u}|^2)\\
&\le 2^{m-2}K(4\delta\int_{\mathbb{B}_1^m(0)}|Du|^2+4\delta^{-1}\int_{\mathbb{B}_1^m(0)} |\xi\circ u-\overline{\xi\circ u}|^2)\\
&=2^mK\delta\int_{\mathbb{B}_1^m(0)}|Du|^2+2^mK\delta^{-1}\int_{\mathbb{B}_1^m(0)}|\xi\circ u-\overline{\xi\circ u}|^2\\
&=\lambda E_1(u)+C\lambda^{-1}\int_{\mathbb{B}_1^m(0)}|\xi\circ u-\overline{\xi\circ u}|^2,
\end{split}
\end{equation*}
where $C=(2^mK)^2$.
\end{proof}
\section{Energy Improvement}
\subsection{A Poincare-Type Theorem}
\begin{definition}
\begin{equation*}
p^*=
\begin{cases} \DF{mp}{m-p} & p<m\\
\mbox{any real number}\in[1,\infty) & p=m\\
\infty & p>m
\end{cases}
\end{equation*}
\end{definition}
\begin{theorem}[$\cite{zw},\;\mbox{Theorem}\; 4.4.6$]
Let $\Omega\subset \mathbb{R}^m$ be a bounded Lipschitz domain and
suppose $u\in W^{1,p}(\Omega)$, $1<p<\infty$. Let
$$c(u)=\int_{\partial \Omega} u d\mathcal{H}^{m-1}$$
Then there is a constant $C=C(m,p,\Omega)$, such that
$$(\int_\Omega |u-c(u)|^{p*}dx)^{1/p^*}\le C(\int_\Omega
|Du|^pdx)^{1/p}.$$
\end{theorem}
\begin{corollary} Let $L$ be a real positive number such that
$\mathcal{H}^{m-1}(\partial \mathbb{B}_L^m(0))=1$. Suppose $u\in
W^{1,2}(\mathbb{B}_L^m(0))$, and let
$$c(u)=\int_{\partial \mathbb{B}_L^m(0)} u d\mathcal{H}^{m-1}=\dashint_{\partial \mathbb{B}_L^m(0)} u d\mathcal{H}^{m-1}$$
Then there is a constant $C=C(m)$ such that
$$\int_{\mathbb{B}_L^m(0)} |u-c(u)|^2dx\le
C\int_{\mathbb{B}_L^m(0)} |Du|^2dx.$$
\end{corollary}
\begin{proof}
Case 1. $m>2$:\\
By H$\ddot{o}$lder inequality, with parameters $m/(m-2)$ and $m/2$
\begin{equation*}
\begin{split}
\int_{\mathbb{B}_L^m(0)}|u-c(u)|^2 dx&\le (\int_{\mathbb{B}_L^m(0)}|u-c(u)|^{2\times \DF{m}{m-2}} dx)^{(m-2)/m}
(\int_{\mathbb{B}_L^m(0)} 1^{\DF{m}{2}} dx)^{2/m}\\
&=(\mathcal{L}^m(\mathbb{B}_L^m(0)))^{2/m}(\int_{\mathbb{B}_L^m(0)}|u-c(u)|^{2^*} dx)^{(m-2)/m}\\
&\le(\mathcal{L}^m(\mathbb{B}_L^m(0)))^{2/m}[C(\int_{\mathbb{B}_L^m(0)}|Du|^2dx)^{1/2}]^{2^*\times (m-2)/m}\\
&=C\int_{\mathbb{B}_L^m(0)}|Du|^2 dx
\end{split}
\end{equation*}
Case 2. $m=2$:\\
We just choose $2^*$ to be $2$ in Theorem 6.1.
\end{proof}
\subsection{Blowing-up Sequence}
\begin{definition}
\begin{equation*}
\begin{split}
\mathcal{F}=\{&u\in\mathcal{Y}_2(\mathbb{B}_1^m(0);\underline{Q}(\mathbb{S}^{n-1})),\forall
0<r<1,
\xi^{-1}\circ \rho\circ AV_{r,0}(\xi\circ u)=Q[[b_r]],\\
&b_r\in\mathbb{R}^n, \DF{1}{2}<|b_r|<\DF{3}{2},\mbox{where}\;
AV_{r,0}(\xi\circ u)=\dashint_{\partial \mathbb{B}_r^m(0)} \xi\circ
u\}
\end{split}
\end{equation*}
\end{definition}
\begin{theorem}[Energy Improvement]
There are positive constants $\epsilon$ and $0<\theta<\DF{1}{2}$
such that if $u$ is a Dir minimizer in $
\mathcal{Y}_2(\mathbb{B}_1^m(0), \underline{Q}(\mathbb{S}^{n-1}))$,
$u\in\mathcal{F}$, $E_1(u)<\epsilon^2$, then $E_\theta(u)\le
\theta^{\omega_{2.13}}E_1(u)$.
\end{theorem}
\begin{proof}
Were the theorem false, there would be, for each $0<\theta<1/2$, a
sequence $u_i\in\mathcal{F}$,$ \epsilon_i^2=E_1(u_i)\to 0$, but
$$E_\theta(u_i)>\theta^{\omega_{2.13}} \epsilon_i^2.$$
Let $\Pi$ be the projection onto the unit sphere in $\mathbb{R}^n$, i.e. $\Pi(x)=\DF{x}{|x|}$. It is easy to check that when we restrict our
attention to the set $U_\epsilon=\{x:1-\epsilon<|x|<1+\epsilon\}$, the Lipschitz constant of $\Pi$ is no more than $1/(1-\epsilon)$.\\
Define
$$\Pi_\sharp\circ \xi^{-1}\circ\rho\circ AV_{L,0}(\xi\circ u_i)=Q[[b_i]],$$
where $L$ is defined in Corollary 6.1.\\ Consider the following
blowing-up sequence
$$\DF{u_i-Q[[b_i]]}{\epsilon_i}.$$
The energy of each one is one by the definition of $\epsilon_i$.
As for their $L^2$ norms, we estimate as follows:
\begin{equation*}
\begin{split}
\mathcal{G}(u_i,Q[[b_i]]) &=|\Pi_\sharp\circ \xi^{-1}\circ \rho\circ\xi\circ u_i-\Pi_\sharp\circ\xi^{-1}\circ\rho\circ\dashint_{\partial
\mathbb{B}_L^m(0)} \xi\circ u_i d\mathcal{H}^{m-1}|\\
&\le
(\mbox{Lip}\Pi|U_{1/2})(\mbox{Lip}\xi^{-1})(\mbox{Lip}\rho)|\xi\circ
u_i-\dashint_{\partial
\mathbb{B}_L^m(0)} \xi\circ u_i d\mathcal{H}^{m-1}|\\
&\le 2(\mbox{Lip}\xi^{-1})(\mbox{Lip}\rho)|\xi\circ u_i-\dashint_{\partial
\mathbb{B}_L^m(0)} \xi\circ u_i d\mathcal{H}^{m-1}|\\
\end{split}
\end{equation*}
From Corollary 6.1, we have
\begin{equation*}
\begin{split}
\int_{\mathbb{B}_L^m(0)}\mathcal{G}(u_i,Q[[b_i]])^2 dx &\le C
\int_{\mathbb{B}_L^m(0)}|\xi\circ u_i-\dashint_{\partial \mathbb{B}_L^m(0)}\xi\circ u_i d\mathcal{H}^{m-1}|^2 dx\\
&\le C\int_{\mathbb{B}_L^m(0)}|D(\xi\circ u_i)|^2 dx\\
&\le C \int_{\mathbb{B}_1^m(0)}|Du_i|^2 dx
\end{split}
\end{equation*}
Hence the $L^2$-norm of the blowing-up sequence in $
\mathbb{B}_L^m(0)$ is uniformly bounded. (Technically, we should
therefore from now on, focus on $\mathbb{B}_L^m(0)$ instead of
$\mathbb{B}_1^m(0)$. But since the regularity is only a local
property, we may just stick to $\mathbb{B}_1^m(0)$ for
convenience.)\\
We use Compactness Theorem 4.2 in $\cite{zw1}$ to get a
subsequence(for convenience, whenever we have to take a
subsequence, we do not change the notation) such that
$$w_i:=\DF{u_i-Q[[b_i]]}{\epsilon_i}\rightharpoonup w\;\mbox{weakly in}\;\mathcal{Y}_2$$
$$w_i\rightarrow w\;\mbox{strongly in}\; L^2$$
$$\int_{\mathbb{B}_1^m(0)}|Dw|^2\le \liminf_{k\to\infty}\int_{\mathbb{B}_1^m(0)}|Dw_i|^2=1,$$
for some $w\in\mathcal{Y}_2(\mathbb{B}_1^m(0),\mathbb{Q})$.
\subsection{Blowing-up the Constraint}
Since $\mathbb{S}^{n-1}$ is compact, we may assume that $b_i\to b\in \mathbb{S}^{n-1}$.\\
Let $u_i(x)=\sum_{j=1}^Q u_j^{(i)}(x)$, $w_i(x)=\sum_{j=1}^Q
w_j^{(i)}(x)$. By definition we have
$$ \DF{u_j^{(i)}-b_i}{\epsilon_i}=w_j^{(i)},$$
hence
$$\DF{u_j^{(i)}}{\epsilon_i}=\DF{b_i}{\epsilon_i}+w_j^{(i)}.$$
Take the norm of both sides,
$$\epsilon_i^{-2}=\epsilon_i^{-2}+|w_j^{(i)}|^2+\DF{2}{\epsilon_i}<b_i,w_j^{(i)}>$$
Hence
$$<b_i,w_j^{(i)}>=-\DF{\epsilon_i}{2}|w_j^{(i)}|^2.$$
Let $i$ go to infinity, we know $w\in
\mathcal{Y}_2(\mathbb{B}_1^m(0),\underline{Q}(P))$ for some $n-1$
dimensional plane $P$ passing through the origin and perpendicular
to $b$.
\subsection{Strong Convergence and Minimality}
Now we want to show that $w$ is Dir minimizing in $\mathcal{Y}_2(
\mathbb{B}_1^m(0),\underline{Q}(P))$ and the convergence of $w_i$
in $\mathcal{Y}_2$ is actually
strong.\\
Let $\mathbb{B}_{\rho_0}^m(y)\subset \mathbb{B}_1^m(0)$, and let
$\delta>0$, $\theta\in(0,1)$ be given. Choose any
$M\in\{1,2,\cdot\cdot\cdot\}$ such that
$$\limsup_{i\to\infty}{\rho_0}^{2-m}\int_{\mathbb{B}_{\rho_0}^m(y)}|D(u_i/\epsilon_i)|^2<M\delta$$
and note that if $\epsilon\in(0,(1-\theta)/M)$, we must have some integer $l\in\{1,2,\cdot\cdot\cdot,M\}$ such that
$$\rho_0^{2-m}\int_{\mathbb{B}_{\rho_0(\theta+l\epsilon)}^m(y)\backslash
\mathbb{B}_{\rho_0(\theta+(l-1)\epsilon)}^m(y)}|D(u_i/\epsilon_i)|^2<\delta
\;\mbox{for infinitely many}\;i$$
This is because that otherwise we get that $\rho_0^{2-m}\int_{\mathbb{B}_{\rho_0}^m(y)}|D(u_i/\epsilon_i)|^2\ge M\delta$ for all sufficiently
large $i$ by summation over $l$, contrary to the definition of $M$. Thus choosing such an $l$, letting $\rho=\rho_0(\theta+(l-1)\epsilon)$, and
noting that $\rho(1+\epsilon)\le \rho_0(\theta+l\epsilon)<\rho_0$, we get $\rho\in[\theta\rho_0,\rho_0)$ such that
$$\rho_0^{2-m}\int_{\mathbb{B}_{\rho(1+\epsilon)}^m(y)\backslash
\mathbb{B}_\rho^m(y)} |D(u_i/\epsilon_i)|^2<\delta\;\mbox{for some
subsequence}\;u_i.$$ By weak convergence, we have
$$\rho_0^{2-m}\int_{\mathbb{B}_{\rho(1+\epsilon)}^m(y)\backslash
\mathbb{B}_\rho^m(y)}
|D(w+\DF{Q[[b_i]]}{\epsilon_i})|^2<\delta\;\mbox{for some
subsequence}.$$ We can not use the Luckhaus-type Theorem 3.2 in
$\cite{zw1}$ now, because $\epsilon_i
w+Q[[b_i]]\not\in\underline{Q}(\mathbb{S}^{n-1})$. But we can use
the technique $``(\Pi_a|\mathbb{S}^{n-1})^{-1}\circ \Pi_a''$ as we
did in proving the hybrid inequality to get a map denoted as
$$(\Pi_i)_\sharp\circ(\epsilon_i w+Q[[b_i]])\in \mathcal{Y}_2(
\mathbb{B}_1^m(0),\underline{Q}(\mathbb{S}^{n-1}))$$
such that
$$\rho_0^{2-m}\int_{
\mathbb{B}_{\rho(1+\epsilon)}^m(y)\backslash
\mathbb{B}_\rho^m(y)}|D((\Pi_i)_\sharp\circ(\epsilon_i
w+Q[[b_i]]))|^2\le \;\mbox{some constant}\cdot \delta$$ Now by
Corollary 3.1(2) in $\cite{zw1}$, since $\int_{
\mathbb{B}_{\rho_0}^m(y)}
\mathcal{G}(u_i,(\Pi_i)_\sharp\circ(\epsilon_i
w+Q[[b_i]]))^2\rightarrow 0$, for sufficiently large $i$, we can
find $v_{i}\in
\mathcal{Y}_2(\mathbb{B}_{\rho(1+\epsilon)}^m(y)\backslash
\mathbb{B}_\rho^m(y);\underline{Q}(\mathbb{S}^{n-1}))$ such that
$v_{i}=(\Pi_i)_\sharp\circ(\epsilon_i w+Q[[b_i]])$ in a
neighborhood of $\partial \mathbb{B}_\rho^m(y)$, $v_{i}=u_{i}$ in
a neighborhood of $\partial \mathbb{B}_{\rho(1+\epsilon)}^m(y)$
and
\begin{equation*}
\begin{split}
\rho^{2-m}\int_{\mathbb{B}_{\rho(1+\epsilon)}^m(y)\backslash \mathbb{B}_\rho^m(y)} |Dv_{i}|^2&\le\\
&C\rho^{2-m}\int_{\mathbb{B}_{\rho(1+\epsilon)}^m(y)\backslash
\mathbb{B}_\rho^m(y)} (|D((\Pi_i)_\sharp\circ(\epsilon_i
w+Q[[b_i]])|^2\\&+ |Du_{i}|^2
+\DF{\mathcal{G}(u_{i},(\Pi_i)_\sharp\circ(\epsilon_i
w+Q[[b_i]]))^2}{\epsilon^2\rho^2}),
\end{split}
\end{equation*}
where $C$ depends only on $m,n,Q$.\\
Now let $v\in
\mathcal{Y}_2(\mathbb{B}_{\theta\rho_0}^m(y),\underline{Q}(P))$
such that $v=w$ in a neighborhood of $\partial
\mathbb{B}_{\theta\rho_0}^m(y)$.\\
Define
$$\tilde{v}=(\Pi_i)_\sharp\circ(\epsilon_i v+Q[[b_i]])\;\mbox{in}\;\mathbb{B}_{\theta\rho_0}^m(y)$$
$$\tilde{v}=(\Pi_i)_\sharp\circ(\epsilon_i w+Q[[b_i]])\;\mbox{in}\;
\mathbb{B}_{\rho_0}^m(y)\backslash
\mathbb{B}_{\theta\rho_0}^m(y).$$ Let $\tilde{u_i}$ be defined by
$$\tilde{u_i}=\tilde{v}\;\mbox{in}\;\mathbb{B}_\rho^m(y),$$
$$\tilde{u_i}=v_i\;\mbox{in}\;\mathbb{B}_{(1+\epsilon)\rho}^m(y)\backslash \mathbb{B}_\rho^m(y),$$
$$\tilde{u_i}=u_i\;\mbox{in}\;\mathbb{B}_{\rho_0}^m(y)\backslash \mathbb{B}_{(1+\epsilon)\rho}^m(y).$$
By the minimizing property of $u_i$, we have
\begin{equation*}
\begin{split}
\int_{\mathbb{B}_{(1+\epsilon)\rho}^m(y)}|Du_i|^2&\le\int_{
\mathbb{B}_{(1+\epsilon)\rho}^m(y)}|D\tilde{u_i}|^2\\
&=\int_{\mathbb{B}_\rho^m(y)}|D\tilde{v}|^2+\int_{
\mathbb{B}_{(1+\epsilon)\rho}^m(y)\backslash \mathbb{B}_\rho^m(y)}
|Dv_i|^2.
\end{split}
\end{equation*}
Therefore
\begin{equation*}
\begin{split}
\rho^{2-m}\int_{\mathbb{B}_\rho^m(y)}|Dw|^2
&\le\liminf_{i\to\infty} \rho^{2-m}\int_{\mathbb{B}_\rho^m(y)}\DF{|Du_i|^2}{\epsilon_i^2}\\
&\le\liminf_{i\to\infty}\rho^{2-m}\int_{\mathbb{B}_\rho^m(y)}\DF{|D\tilde{v}|^2}{\epsilon_i^2}+\liminf_{i\to\infty}
\rho^{2-m}\int_{\mathbb{B}_{(1+\epsilon)\rho}^m(y)\backslash \mathbb{B}_\rho^m(y)} \DF{|Dv_i|^2}{\epsilon_i^2}\\
&\le \rho^{2-m}\int_{\mathbb{B}_{\theta\rho_0}^m(y)} |Dv|^2+\rho^{2-m}\int_{\mathbb{B}_\rho^m(y)\backslash
\mathbb{B}_{\theta\rho_0}^m(y)} |Dw|^2+C\delta
\end{split}
\end{equation*}
Since $\delta$ was arbitrary, we have
$$\rho^{2-m}\int_{\mathbb{B}_{\theta\rho_0}(y)}|Dw|^2\le \rho^{2-m}\int_{\mathbb{B}_{\theta\rho_0}(y)}|Dv|^2$$
Therefore, $w$ is minimizing on $\mathbb{B}_{\theta\rho_0}^m(y)$,
and in view of the arbitrariness of $\theta$ and $\rho_0$, this
shows that $w$ is minimizing on all balls $\mathbb{B}_\rho^m(y)$
with $
\mathbb{B}_\rho^m(y)\subset \mathbb{B}_1^m(0)$.\\
Finally to prove that the convergence is strong we note that if we use $v=w$ as above, we can conclude
$$\liminf_{i\to\infty}\rho^{2-m}\int_{\mathbb{B}_\rho^m(y)}\DF{|Du_i|^2}{\epsilon_i^2}\le \rho^{2-m}\int_{
\mathbb{B}_\rho^m(y)}|Dw|^2+C\delta$$
and hence, in view of the arbitrariness of $\theta$ and $\delta$,
$$\rho^{2-m}\liminf_{i\to\infty}\int_{\mathbb{B}_{\rho_1}^m(y)}\DF{|Du_i|^2}{\epsilon_i^2}\le \rho^{2-m}\int_{
\mathbb{B}_{\rho_0}^m(y)} |Dw|^2,$$
for each $\rho_1<\rho_0$. Evidently it follows from this(keeping in mind the arbitrariness of $\rho_0$) that
$$\liminf_{i\to\infty}\int_{\mathbb{B}_\rho^m(y)}\DF{|Du_i|^2}{\epsilon_i^2}\le \int_{\mathbb{B}_\rho^m(y)}|Dw|^2$$
for every ball $\mathbb{B}_\rho^m(y)$ such that $
\mathbb{B}_\rho^m(y)\subset \mathbb{B}_1^m(0)$. Then since
$$\int_{\mathbb{B}_\rho^m(y)}|D(u_i/\epsilon_i)-Dw|^2=\int_{
\mathbb{B}_\rho^m(y)}|Dw|^2+\int_{\mathbb{B}_\rho^m(y)}|D(u_i/\epsilon_i)|^2-2
\int_{\mathbb{B}_\rho^m(y)}Dw\cdot D(u_i/\epsilon_i),$$ we can
evidently select a subsequence which converges strongly to $Dw$ on
$\mathbb{B}_\rho^m(y)$. Since this holds for arbitrary
$\mathbb{B}_\rho^m(y)\subset \mathbb{B}_1^m(0)$, it is then easy
to see(by covering $\mathbb{B}_1^m(0)$ by a countable collection
of balls $\mathbb{B}_{\rho_j}^m(y_j)$ with $
\mathbb{B}_{\rho_j}^m(y_j)\subset \mathbb{B}_1^m(0)$) that there
is a subsequence such that $D(u_i/\epsilon_i)$ converges strongly
locally in all of $\mathbb{B}_1^m(0)$.
\subsection{Proof of Energy Improvement}
Let's estimate $\dashint_{\mathbb{B}_r^m(0)}|\xi\circ u_i-\overline{\xi\circ u_i}|^2dx$,where
$\overline{\xi\circ u_i}=\dashint_{\mathbb{B}_r^m(0)} \xi\circ u_i dx$.
\begin{equation*}
\begin{split}
\dashint_{\mathbb{B}_r^m(0)}|\xi\circ u_i-\overline{\xi\circ u_i}|^2 dx &\le Cr^{2-m}\int_{
\mathbb{B}_r^m(0)}|D(\xi\circ u_i)|^2 dx\;(\mbox{by Poincare
inequality})\\
&=Cr^{2-m}\int_{\mathbb{B}_r^m(0)}|Du_i|^2 dx\\&=Cr^{2-m}\epsilon_i^2\int_{\mathbb{B}_r^m(0)}|D(u_i/\epsilon_i)|^2 dx\\
&=Cr^{2-m}\epsilon_i^2\int_{\mathbb{B}_r^m(0)} |Dw_i|^2 dx
\end{split}
\end{equation*}
We have already proved that $Dw_i$ converges strongly to $Dw$ in
$\mathcal{Y}_2$, hence
$$\dashint_{\mathbb{B}_r^m(0)}|\xi\circ u_i-\overline{\xi\circ u_i}|^2dx\le Cr^{2-m}\epsilon_i^2\int_{
\mathbb{B}_r^m(0)}|Dw|^2 dx.$$ We also have proved the Dir
minimality of $w$, hence by Theorem 2.4
$$\dashint_{\mathbb{B}_r^m(0)}|\xi\circ u_i-\overline{\xi\circ u_i}|^2dx\le Cr^{2-m}\epsilon_i^2r^{m-2+2\omega_{2.13}}
=Cr^{2\omega_{2.13}}\epsilon_i^2.$$ Applying the Hybrid Inequality
to $u_i(2\theta x)$, we get
$$E_{1/2}(u_i(2\theta x))\le \lambda E_1(u_i(2\theta x))+C\lambda^{-1}
\int_{\mathbb{B}_1^m(0)}|\xi\circ u_i(2\theta
x)-\overline{\xi\circ u_i(2\theta x)}|^2 dx$$ which can be
simplified to
$$E_\theta(u_i)\le \lambda E_{2\theta}(u_i)+C\lambda^{-1}\dashint_{\mathbb{B}_{2\theta}^m(0)}|\xi
\circ u_i-\overline{\xi\circ u_i}|^2 dx$$
$$\le \lambda E_{2\theta}(u_i)+C\lambda^{-1}\cdot(2\theta)^{2\omega_{2.13}}\epsilon_i^2$$
Choosing the positive integer $k=k(\theta)$ for which $1/2\le 2^k\theta\le 1$, we iterate $k-1$ more times to obtain
\begin{equation*}
\begin{split}
E_\theta(u_i)&\le \lambda^k E_{2^k\theta}(u_i)+\sum_{j=1}^k \lambda^{j-1}C\lambda^{-1}\dashint_{
\mathbb{B}_{2^j\theta}^m(0)}|
\xi\circ u_i-\overline{\xi\circ u_i}|^2\\
&\le \lambda^k(1/2)^{2-m}\epsilon_i^2+\sum_{j=1}^k \lambda^{j-1}C\lambda^{-1}C(2^j\theta)^{2\omega_{2.13}}\epsilon_i^2\\
&\le \lambda^k\cdot 2^{m-2}\epsilon_i^2+\sum_{j=1}^\infty (\lambda\cdot 2^{2\omega_{2.13}})^j C\lambda^{-2}\theta^{2\omega_{2.13}}\epsilon_i^2\\
&\le [\lambda^k\cdot 2^{m-2}+\DF{\lambda\cdot 2^{2\omega_{2.13}}}{1-\lambda\cdot 2^{2\omega_{2.13}}}C\lambda^{-2}\theta^{2\omega_{2.13}}]\epsilon_i^2
\end{split}
\end{equation*}
Take $\lambda=\theta^{\DF{m+\omega_{2.13}}{k}}$, we have $
\lambda^k\cdot 2^{m-2}=\theta^{m+\omega_{2.13}}\cdot 2^{m-2}=\theta^m\cdot 2^{m-2}\cdot \theta^{\omega_{2.13}}\le
(1/2)^m\cdot 2^{m-2}\theta^{\omega_{2.13}}\le \theta^{\omega_{2.13}}/4.$\\
Since $\lambda=\theta^{\DF{m+\omega_{2.13}}{k}}\le (2^{-k})^{\DF{m+\omega_{2.13}}{k}}=2^{-(m+\omega_{2.13})}$,
\begin{equation*}
\begin{split}
\DF{\lambda\cdot 2^{2\omega_{2.13}}}{1-\lambda\cdot 2^{2\omega_{2.13}}}C\lambda^{-2}\theta^{2\omega_{2.13}}
&\le \DF{2^{2\omega_{2.13}}C}{1-2^{\omega_{2.13}-m}}\theta^{-\DF{m+\omega_{2.13}}{k}}\theta^{2\omega_{2.13}}\\
&=
K\theta^{\omega_{2.13}-\DF{m+\omega_{2.13}}{k}}\theta^{\omega_{2.13}}
\end{split}
\end{equation*}
Let's choose $\theta$ small enough such that $\theta^{\omega_{2.13}-\DF{m+\omega_{2.13}}{k}}\le 1/4K$.
This is possible because it is equivalent to
$$\theta^{\omega_{2.13}}\le \theta^{\DF{m+\omega_{2.13}}{k}}/4K$$
Noting that $\theta\ge 2^{-1-k}$, the right side of above one is greater than
$$2^{-(k+1)(m+\omega_{2.13})/k}/4K$$
which is bounded from below although when $\theta$ goes to zero, $k$ goes to infinity.\\
Thus for $i$ sufficiently large enough,we have
$$E_\theta(u_i)\le(\DF{1}{4}\theta^{\omega_{2.13}}+\DF{1}{4}\theta^{\omega_{2.13}})\epsilon_i^2
<\theta^{\omega_{2.13}}\epsilon_i^2,$$
contradicting the choice of $u_i$.
\end{proof}
\section{Energy Decay}
\begin{theorem}[Energy decay] If $u\in \mathcal{F}$ is Dir minimizing, $\mathbb{B}_R^m(0)\subset
\mathbb{B}_1^m(0)$, and $R^{2-m}\int_{\mathbb{B}_R^m(0)}|Du|^2\le
\epsilon^2$, then
$$\int_{\mathbb{B}_r^m(0)}|Du|^2\le\theta^{2-m-\omega_{2.13}}R^{-\omega_{2.13}}\epsilon^2r^{m-2+\omega_{2.13}},\;\mbox{for}\;0\le r\le R$$
where $\epsilon$ and $\theta$ are as in the Energy Improvement.
\end{theorem}
\begin{proof}Let
$u_{\theta^iR}\equiv u(\theta^i Rx),i=0,1,\cdot\cdot\cdot$. It is easy to check
$$E_1(u_{\theta^k R})=(\theta^k R)^{2-m}Dir(u,\mathbb{B}_{\theta^k R}^m(0))=E_\theta(u_{\theta^{k-1}R})$$
Claim:
$E_\theta(u_R)\le \theta^{\omega_{2.13}} \epsilon^2$.\\
This is because $u_R\in\mathcal{F}$,and
$E_1(u_R)=R^{2-m}Dir(u,\mathbb{B}_R^m(0))\le \epsilon^2$ by our
assumption. Hence we can use the
energy improvement to the function $u_R$ to get that.\\
Claim: $E_\theta(u_{\theta R})\le
\theta^{2\omega_{2.13}}\epsilon^2$.\\ Obviously, $u_{\theta
R}\in\mathcal{F}$, moreover,
$$E_1(u_{\theta R})=E_\theta(u_R)\le \theta^{\omega_{2.13}} E_1(u_R)\le \theta^{\omega_{2.13}}\epsilon^2\le \epsilon^2.$$
Hence using the energy improvement to function $u_{\theta R}$, we
get
$$E_\theta(u_{\theta R})\le \theta^{\omega_{2.13}} E_1(u_{\theta R})\le \theta^{2\omega_{2.13}}\epsilon^2.$$
Continuing the process, we get
$$E_1(u_{\theta^k R})=E_\theta(u_{\theta^{k-1}R})\le
\theta^{\omega_{2.13}} E_1(u_{\theta^{k-1}R})
=\theta^{\omega_{2.13}} E_\theta(u_{\theta^{k-2}R})$$
$$\le \theta^{2\omega_{2.13}}E_1(u_{\theta^{k-2}R})\cdot\cdot\cdot=\theta^{k\omega_{2.13}}\epsilon^2,$$
for $k=1,2,3,\cdot\cdot\cdot.$ \\
Given $0<r\le R$, choose $k$ such that $\theta^{k+1}R<r\le
\theta^k R$.
\begin{equation*}
\begin{split}
r^{2-m}\int_{\mathbb{B}_r^m(0)}|Du|^2&\le (\theta^{k+1}R)^{2-m}\int_{
\mathbb{B}_{\theta^k R}^m(0)}|Du|^2\\&=
\theta^{2-m}(\theta^k R)^{2-m}\int_{\mathbb{B}_{\theta^k R}^m(0)} |Du|^2\\
&=\theta^{2-m}E_1(u_{\theta^k R})\\
&\le \theta^{2-m}\theta^{k\omega_{2.13}}\epsilon^2\\&=\theta^{2-m-\omega_{2.13}}
\theta^{(k+1)\omega_{2.13}}\epsilon^2\\
&\le\theta^{2-m-\omega_{2.13}}(r/R)^{\omega_{2.13}}\epsilon^2
\end{split}
\end{equation*}
\end{proof}
\section{$\mathcal{H}^{m-2}(\mbox{singular set})=0$}
\begin{theorem}
Let
$u\in\mathcal{Y}_2(\mathbb{B}_1^m(0),\underline{Q}(\mathbb{S}^{n-1}))$
be a strictly defined, Dirichlet minimizing map. Then it is
H$\ddot{o}$lder continuous away from the boundary except for a
closed subset $S\subset \mathbb{B}_1^m(0)$ such that
$\mathcal{H}^{m-2}(S)=0$.
\end{theorem}
\begin{proof}
Let $$S=\{x\in\mathbb{B}_1^m(0),\limsup_{\rho\downarrow 0}
\rho^{2-m}\int_{\mathbb{B}_\rho^m(x)} |Du|^2>0\}.$$ Obviously, $S$
is closed and $\mathcal{H}^{m-2}(S)=0$ (see for example Lemma
2.1.1 in $\cite{ly}$). \\Let's look at a point $a\in
\mathbb{B}_1^m(0)\sim S$. We may assume
$a=0$.\\
Let $\epsilon$ be the constant in the Energy Improvement, and
$k=k(Q,m,n)$ be the constant in the ``small energy regularity"
theorem in $\cite{lc1}$.  Since $0\notin S$, there is $R>0$ such
that $\mathbb{B}_{2R}^m(0)\subset \mathbb{B}_1^m(0)\sim S$ and
$$R^{2-m}\int_{\mathbb{B}_{2R}^m(0)}|Du|^2\le \min\{\epsilon^2,k\}.$$
For any $b\in \mathbb{B}_R^m(0)$,
$$R^{2-m}\int_{\mathbb{B}_R^m(b)}|Du|^2\le R^{2-m}\int_{\mathbb{B}_{2R}^m(0)}|Du|^2\le \min\{\epsilon^2,k\}.$$
We have two possibilities:\\
Case 1: $b\not\in B_0$. By the ``small energy regularity" theorem
in $\cite{lc1}$, we have
$$\int_{\mathbb{B}_r^m(b)}|Du|^2\le \mbox{some constant}\cdot r^{m-2+\beta},0\le r\le R$$
where $\beta$ is the constant given in $\cite{lc1}$.\\
Case 2: $b\in B_0$. From Energy Decay we have
$$\int_{\mathbb{B}_r^m(b)}|Du|^2\le \theta^{2-m-\omega_{2.13}}R^{-\omega_{2.13}}\epsilon^2r^{m-2+\omega_{2.13}},0\le r\le R.$$
Therefore $u\in
C^{0,\min\{\omega_{2.13},\beta\}/2}[\mathbb{B}_R^m(0)]$ by
Morrey's growth lemma.
\end{proof}
\section{Dimension Reduction}
\subsection{Monotonicity Formula}
Suppose $u:\mathbb{B}_1^m(0)\to \underline{Q}(\mathbb{S}^{n-1})\subset\underline{Q}(\mathbb{R}^n)$ is Dir minimizing, although $\xi\circ u$ is not necessarily harmonic,
we still have the following results:\\
Consider the domain variation
$$u_s(x)=u(x+s\zeta(x)), \;\mbox{where}\;\zeta=(\zeta^1,\cdot\cdot\cdot,\zeta^m),\;\mbox{with}\;\zeta^j\in C_c^\infty(
\mathbb{B}_1^m(0)).$$
We should have
$$\DF{d}{ds}|_{s=0}\;\mbox{Energy of}\;\xi\circ u_s=0.$$
If we let $f=\xi\circ u$, it is easy to check as in $\cite{sl},\S 2.2$
$$\int_{\mathbb{B}_1^m(0)}\sum_{i,j=1}^m(|Df|^2\delta_{ij}-2D_i f\cdot D_j f)D_i \zeta^j=0.$$
\begin{theorem}[Monotonicity Formula]
If $u:\mathbb{B}_1^m(0)\to \underline{Q}(\mathbb{S}^{n-1})\subset \underline{Q}(\mathbb{R}^n)$ is Dir minimizing, then
$$\rho^{2-m}\int_{\mathbb{B}_\rho^m(x)}|Du|^2-\sigma^{2-m}\int_{\mathbb{B}_\sigma^m(x)}|Du|^2=2\int_{
\mathbb{B}_\rho^m(x)\backslash \mathbb{B}_\sigma^m(x)}R^{2-m}
|\DF{\partial u}{\partial R}|^2$$ for any $0<\sigma<\rho<\rho_0$,
provided $\mathbb{B}_{\rho_0}^m(x)\subset \mathbb{B}_1^m(0)$,
where $R=|y-x|$ and $
\partial/\partial R$ means directional derivative in the radial direction $|y-x|^{-1}(y-x)$.
\end{theorem}
\begin{proof}
Just apply the argument of ($\cite{sl},\S 2.4$) to the function $\xi\circ u$ to get
$$\rho^{2-m}\int_{\mathbb{B}_\rho^m(x)}|D(\xi\circ u)|^2-\sigma^{2-m}\int_{\mathbb{B}_\sigma^m(x)}|D(\xi\circ u)|^2
=2\int_{\mathbb{B}_\rho^m(x)\backslash
\mathbb{B}_\sigma^m(x)}R^{2-m} |\DF{\partial (\xi\circ
u)}{\partial R}|^2$$ and notice that $|D_v(\xi\circ u)|=|D_v u|$.
\end{proof}
\begin{remark} (1) From above, $\rho^{2-m}\int_{\mathbb{B}_\rho^m(x)} |Du|^2$
is an increasing function of $\rho$ for $\rho\in (0,\rho_0)$, and
hence the limit as $\rho\to 0$ of
$\rho^{2-m}\int_{\mathbb{B}_\rho^m(x)}|Du|^2$ exists; this limit
is denoted as $\Theta_u(x)$.
It is also easy to see that the density $\Theta_u$ is upper semi-continuous on $\mathbb{B}_1^m(0)$.\\
(2) Another important additional conclusion, which we see by
taking the limit as $\sigma\to 0$ in the monotonicity formula, is
that $\int_{\mathbb{B}_\rho^m(x)}R^{2-m}|\DF{\partial u}{\partial
R}|^2<\infty$ and
$$\rho^{2-m}\int_{\mathbb{B}_\rho^m(x)}|Du|^2-\Theta_u(x)=2\int_{
\mathbb{B}_\rho^m(x)}R^{2-m}|\DF{\partial u}{\partial R}|^2.$$
\end{remark}
\subsection{Definition of Tangent Maps}
Let $\mathbb{B}_{\rho_0}^m(y)$ with $
\mathbb{B}_{\rho_0}^m(y)\subset \mathbb{B}_1^m(0)$, and for any
$\rho>0$ consider the scaled function $ u_{y,\rho}$ defined by
$$u_{y,\rho}(x)=u(y+\rho x).$$
If $\sigma>0$ is arbitrary and $\rho<\rho_0/\sigma$, we have
(using $Du_{y,\rho}(x)=\rho(Du)(y+\rho x)$, and making a change of
variable $\tilde{x}=y+\rho x$ in the energy integral of
$u_{y,\rho}$)
\begin{equation}
\sigma^{2-m}\int_{\mathbb{B}_\sigma^m(0)}|Du_{y,\rho}|^2=(\sigma\rho)^{2-m}\int_{
\mathbb{B}_{\sigma\rho}^m(y)} |Du|^2\le
{\rho_0}^{2-m}\int_{\mathbb{B}_{\rho_0}^m(y)}|Du|^2
\end{equation}
Thus if $\rho_j\downarrow 0$ then $\limsup_{j\to\infty}\int_{\mathbb{B}_\sigma^m(0)}|Du_{y,\rho_j}|^2<\infty$ for each $\sigma>0$.
Their $L^2-$norms
$$\int_{\mathbb{B}_\sigma^m(0)}|u_{y,\rho}|^2=\rho^{-m}\int_{\mathbb{B}_{\sigma\rho}^m(y)}|u|^2<\infty$$
uniformly for $\rho$ because $u(x)\in \underline{Q}(\mathbb{S}^{n-1})$.\\
So we can use Compactness Theorem 4.3 in $\cite{zw1}$ to get a
subsequence $\rho_{j'}$ such that $u_{y,\rho_{j'}}\to \varphi$
locally in $\mathbb{R}^m$ with respect to the
$\mathcal{Y}_2-$norm, where $\varphi:\mathbb{R}^m\to
\underline{Q}(\mathbb{S}^{n-1})$ is an energy minimizing map, called a tangent map of $u$ at $y$.\\
\subsection{Properties of Tangent Maps}
Let $\rho_j\downarrow 0$ be one of the sequences such that the
re-scaled maps $u_{y,\rho_j}\to \varphi$ as described above. Since
$u_{y,\rho_j}$ converges in energy to $\varphi$, we have, after
setting $\rho=\rho_j$ and taking limits on each side of (1) as
$j\rightarrow\infty$,
$$\sigma^{2-m}\int_{\mathbb{B}_\sigma^m(0)}|D\varphi|^2=\Theta_u(y).$$
Thus in particular, $\sigma^{2-m}\int_{\mathbb{B}_\sigma^m(0)}|D\varphi|^2$ is a constant function of $\sigma$ and since by definition
$\Theta_\varphi(0)=\lim_{\sigma\downarrow 0}\sigma^{2-m}\int_{\mathbb{B}_\sigma^m(0)} |D\varphi|^2$, we have
\begin{equation}
\Theta_u(y)=\Theta_\varphi(0)\equiv\sigma^{2-m}\int_{\mathbb{B}_\sigma^m(0)}|D\varphi|^2,\forall
\sigma>0
\end{equation}
Thus any tangent map of $u$ at $y$ has constant scaled energy and equal to the density of $u$ at $y$.\\
Furthermore, we apply the monotonicity formula to $\varphi$ to get
$$0=\sigma^{2-m}\int_{\mathbb{B}_\sigma^m(0)}|D\varphi|^2-\tau^{2-m}\int_{
\mathbb{B}_\tau^m(0)}|D\varphi|^2=
\int_{\mathbb{B}_\sigma^m(0)\backslash
\mathbb{B}_\tau^m(0)}R^{2-m}|\DF{\partial \varphi}{\partial
R}|^2.$$ So that $\partial \varphi/\partial R=0$ a.e, and since
$\varphi\in
\mathcal{Y}_2(\mathbb{R}^m,\underline{Q}(\mathbb{S}^{n-1}))$ it is
correct to conclude from this, by integration along rays, that
$$\varphi(\lambda x)\equiv\varphi(x)\;\forall\lambda>0,x\in\mathbb{R}^m$$
\begin{theorem} $y\in\;\mbox{reg}\;u\Leftrightarrow
\Theta_u(y)=0\Leftrightarrow \exists \;\mbox{a constant tangent
map}\;\varphi\; \mbox{of}\;u\;\mbox{at}\;y$
\end{theorem}
\begin{proof}
The first part of the statement is easily obtained from Theorem
8.1. The second part comes from (2).
\end{proof}
\subsection{Properties of Homogeneous Degree Zero Minimizers}
Suppose $\varphi:\mathbb{R}^m\to\underline{Q}(\mathbb{S}^{n-1})$ is a homogeneous degree zero minimizer. We first observe that the density $
\Theta_\varphi(y)$ is maximum at $y=0$; in fact, by the monotonicity formula, for each $\rho>0$ and each $y\in\mathbb{R}^m$
$$2\int_{\mathbb{B}_\rho^m(y)}R_y^{2-m}|\DF{\partial \varphi}{\partial R_y}|^2+\Theta_\varphi(y)=
\rho^{2-m}\int_{\mathbb{B}_\rho^m(y)}|D\varphi|^2,$$
where $R_y(x)\equiv|x-y|$ and $\partial/\partial R_y=|x-y|^{-1}(x-y)\cdot D$. Now $\mathbb{B}_\rho^m(y)^\subset
\mathbb{B}_{\rho+|y|}^m(0)$, so that
\begin{equation*}
\begin{split}
\rho^{2-m}\int_{\mathbb{B}_\rho^m(y)}|D\varphi|^2 &\le \rho^{2-m}\int_{\mathbb{B}_{\rho+|y|}^m(0)}|D\varphi|^2\\
&=(1+\DF{|y|}{\rho})^{m-2}(\rho+|y|)^{2-m}\int_{\mathbb{B}_{\rho+|y|}^m(0)} |D\varphi|^2\\
&\equiv (1+\DF{|y|}{\rho})^{m-2}\Theta_\varphi(0)
\end{split}
\end{equation*}
Thus letting $\rho\uparrow \infty$, we get
$$2\int_{\mathbb{R}^m}R_y^{2-m}|\DF{\partial \varphi}{\partial R_y}|^2+\Theta_\varphi(y)\le \Theta_\varphi(0),$$
which establishes the required inequality
$$\Theta_\varphi(y)\le \Theta_\varphi(0).$$
Notice also that this argument shows that the equality implies that $\partial \varphi/\partial R_y=0$ a.e; that is
$\varphi(y+\lambda x)\equiv\varphi(y+x)$ for each $\lambda>0$. Since we also have $\varphi(\lambda x)\equiv\varphi(x)$ we can then
compute that for any $\lambda>0$ and $x\in \mathbb{R}^m$ that
\begin{equation*}
\begin{split}
\varphi(x)&=\varphi(\lambda x)=\varphi(y+(\lambda x-y))=\varphi(y+\lambda^{-2}(\lambda x-y))\\
&=\varphi(\lambda(y+\lambda^{-2}(\lambda x-y)))=\varphi(x+ty),
\end{split}
\end{equation*}
where $t=\lambda-\lambda^{-1}$ is an arbitrary real number. So let $S(\varphi)$ be defined by
$$S(\varphi)=\{y\in \mathbb{R}^m:\Theta_\varphi(y)=\Theta_\varphi(0)\}.$$
Then we have shown that $\varphi(x)\equiv\varphi(x+ty)$ for all
$x\in \mathbb{R}^m$,$t\in \mathbb{R}$, and $y\in S(\varphi)$.Then
of course $\varphi(x+az_1+bz_2)\equiv \varphi(x)$ for all $a,b\in
\mathbb{R}$ and $z_1,z_2\in S(\varphi)$. But if $z\in
\mathbb{R}^m$ and $ \varphi(x+z)\equiv\varphi(x)$ for all $x\in
\mathbb{R}^m$, then trivially,
$\Theta_\varphi(z)=\Theta_\varphi(0)$(and hence $z\in S(\varphi)$
by definition of $S(\varphi)$), so we conclude
$$S(\varphi)\;\mbox{is a linear subspace of}\;\mathbb{R}^m\;\mbox{and}\;\varphi(x+y)\equiv\varphi(x),x\in\mathbb{R}^m,y\in S(\varphi).$$
(Thus $\varphi$ is invariant under the composition with translation by elements of $S(\varphi)$.) Notice of course that
$$\dim S(\varphi)=n\Leftrightarrow S(\varphi)=\mathbb{R}^m\Leftrightarrow \varphi=\;\mbox{const.}$$
Also, a homogeneous degree zero map which is not constant clearly can not be continuous at $0$, so we always have $0\in
\;\mbox{sing}\;\varphi$ if $\varphi$ is non-constant, and hence, since $\varphi(x+z)\equiv \varphi(x)$ for any $z\in
S(\varphi)$, we have
$$S(\varphi)\subset \;\mbox{sing}\varphi$$
for any non-constant homogeneous degree zero minimizer $\varphi$.
\subsection{Further Properties of sing $u$}
We know
\begin{equation}
y\in\;\mbox{sing}\;u\Leftrightarrow \dim S(\varphi)\le
n-1\;\mbox{for every tangent
map}\;\varphi\;\mbox{of}\;u\;\mbox{at}\;y
\end{equation} Now for each $j=0,1,\cdot\cdot\cdot,n-1$ we define
$$S_j=\{y\in\;\mbox{sing}\;u:\dim S(\varphi)\le j\;\mbox{for all tangent maps}\;\varphi\;\mbox{of}\;u\;\mbox{at}\;y\}$$
Then we have
$$S_0\subset S_1\subset\cdot\cdot\cdot\subset S_{m-3}=S_{m-2}=S_{m-1}=\;\mbox{sing}\;u.$$
To see this first note that $S_{j-1}\subset S_j$ is true by
definition and $S_{m-1}=\;\mbox{sing}\;u$ is just (3). Also, if
$S_{m-3}$ is not equal to both $S_{m-2}$ and $S_{m-1}$, then we
can find $y\in\mbox{sing}\;u$ at which there is a tangent map
$\varphi$ with $\dim S(\varphi)=m-1$ or $m-2$; but then
$\mathcal{H}^{m-2}(S(\varphi))=\infty$ and hence (since
$S(\varphi)\subset \;\mbox{sing}\;\varphi$) we have
$\mathcal{H}^{m-2}(\;\mbox{sing}\;\varphi)=\infty$, contradicting
the fact that $\mathcal{H}^{m-2}(\mbox{sing}\;\varphi)=0$ by
Theorem 8.1. \begin{lemma}For each $j=0,1,\cdot\cdot\cdot,m-3,\dim
S_j\le j$, and, for each $\alpha\ge 0$, $S_0\cap \{
x:\Theta_u(x)=\alpha\}$ is a discrete set. \end{lemma}
\begin{proof}
The proof is exactly the same as in $\cite{sl},\S 3.4$
\end{proof}
\begin{corollary}
Let
$u\in\mathcal{Y}_2(\mathbb{B}_1^m(0),\underline{Q}(\mathbb{S}^{n-1}))$
be a strictly defined, Dirichlet minimizing map. Then it is
H$\ddot{o}$lder continuous away from the boundary except for a
closed subset $S\subset \mathbb{B}_1^m(0)$ such that $\dim(S)\le
m-3$.
\end{corollary}
\begin{proof}
Combine Lemma 9.1 with the fact that sing$(u)=S_{m-3}$.
\end{proof}

\end{document}